\documentclass[12pt]{amsart}

\usepackage{amssymb}
\usepackage{amsmath}
\usepackage{amsfonts}
\usepackage[mathscr]{eucal}
\usepackage{graphicx}
\usepackage{tikz}
\usepackage{hyperref}
\usetikzlibrary{arrows,chains,matrix,positioning,scopes}

\newcommand{\p}{\vskip .5cm}
\newcommand{\pp}{\vskip .17cm}

\newcommand{\s}{\;\;\;\;\;\;\;\;}

\newcommand{\ra}{\rightarrow}

\newcommand{\hra}{\hookrightarrow}
\newcommand{\xra}{\xrightarrow}

\newcommand{\LG}{\mathfrak{g}}
\newcommand{\G}{\mathbf{G}}
\newcommand{\Z}{\mathbb{Z}}
\newcommand{\W}{\mathcal{W}}
\newcommand{\CC}{\mathcal{C}}
\newcommand{\CP}{\mathcal{P}}

\newcommand{\CV}{\mathcal{V}}
\newcommand{\V}{\mathcal{V}}
\newcommand{\CO}{\mathcal{O}}
\newcommand{\sG}{\mathsf{G}}
\newcommand{\sU}{\mathsf{U}}
\newcommand{\se}{\mathsf{e}}
\newcommand{\Lie}{\text{Lie }}

\newcommand{\Ad}{\text{Ad}}

\newcommand{\bsl}{\backslash}
\newcommand{\dsp}{\displaystyle}
\newcommand{\til}{\tilde}

\newtheorem{theorem}{Theorem}[section]
\newtheorem{lemma}[theorem]{Lemma}
\newtheorem{claim}[theorem]{Claim}

\newtheorem{definition}[theorem]{Definition}
\newtheorem{hypothesis}[theorem]{Hypothesis}
\theoremstyle{remark}
\newtheorem{remark}[theorem]{Remark}

\numberwithin{equation}{section}

\begin{document}
\title{Computations of orbital integrals and Shalika germs}
\author{Cheng-Chiang Tsai}
\begin{abstract}
For a reductive group $G$ over a non-archimedean local field, with some assumptions on (residue) characteristic we give an method to compute certain orbital integrals using a method close to that of \cite{GKM06} but in a different language. These orbital integrals allow us to compute the Shalika germs at some ``very elliptic'' elements in terms of number of rational points on some quasi-finite covers of the Hessenberg varieties in \cite{GKM06} which are subvarieties of (partial) flag varieties. Such values of Shalika germs determine the Harish-Chandra local character expansions of the so-called very supercuspidal representations.
\end{abstract}
\maketitle
\vskip -0.5cm

\tableofcontents

\section{Introduction}
\label{intro}

Throughout the article, let $F$ be a non-archimedean local field with residue field $k$ of order $q$. Let $\G$ be a connected reductive group over $F$. We'll actually work with the assumption that $\G$ is semi-simple and simply connected (this makes no harm; see beginning of Sec. \ref{orb}). The $F$-points of $\G$ will be denoted by $G=\G(F)$, and the $F$-points of the Lie algebra of $G$ by $\LG$. Similar conventions will be applied for other algebraic groups, which will be defined either over $F$ or $k$. We will assume that $\text{char}(k)$ is very good for $\G$, $\G$ is tamely ramified and $\text{char}(F)$ is sufficiently large (see Appendix).\p

For any element $g\in\G$, let $C_G(g)$ denote the centralizer of $g$ in $G$ and $\Ad(G)(g)\cong G/C_G(g)$ the orbit of $g$. $\Ad(G)(g)$ is a $p$-adic manifold equipped with a $G$-invariant measure \cite[III.3.27]{SS70}. One may consider the orbital integral $\mu_g(f)=\int_{\Ad(G)(g)}f(x)dx=\int_{G/C_G(g)}f(\Ad(g)x)$. Orbital integrals on the Lie algebra are defined in the same manner. For regular semisimple $g$, such orbital integrals appear on the geometric side of the trace formula and are fundamental objects of study in the Langlands program.\p

For representations of $G$ we have the local character expansion (By Howe \cite{Ho} for $GL_n$ over a $p$-adic field, Harish-Chandra \cite{HC} for $p$-adic fields and Adler-Korman \cite{AK07} in general) for representations of $G$. Let $\pi$ be any irreducible admissible representation of $G$ and $\gamma\in G$ be semisimple. Let $\mathbf{M}=(C_{\G}(\gamma))^o$ be its connected centralizer, $M=\mathbf{M}(F)$ and $\mathfrak{m}=\Lie M$. Then there exists a neighborhood $U$ of $0\in\mathfrak{m}$ such that
\begin{equation}\label{HCHlocal}\Theta_{\pi}(\gamma\cdot\se(X))=\sum_{\CO\in\CO(\gamma)}c_{\CO}(\pi)\hat{\mu}_{\CO}(X),\;X\in U
\end{equation}

where $\Theta_{\pi}$ is the character of $\pi$, $\se$ is a suitable generalization of the exponential map (or just the exponential map when $\text{char}(F)=0$), and $\CO(\gamma)$ is the set of nilpotent orbits in $\mathfrak{m}$, with $c_{\CO}\in\mathbb{C}$ and $\hat{\mu}_{\CO}$ the Fourier transform of $\mu_{\CO}$. Both side of the equality may be viewed either as distributions on $U$, or as locally constant functions defined almost everywhere that represent the distributions.\p

In fact, a general philosophy dated back to Harish-Chandra is that characters should resemble Fourier transforms of orbital integral. Another important result of this philosophy is the work of Kim and Murnaghan \cite{KM06}. Roughly speaking, Kim and Murnaghan showed that whenever an irreducible admissible representation has certain $K$-type, the local character near the identity can be written as a linear combination of Fourier transforms of orbital integrals (not necessarily semisimple or nilpotent). Adler and Spice \cite{AS09} generalize their results in the case of {\it very supercuspidal} representations (see below) to give full a character formula, where the character is locally given as Fourier transforms of elliptic orbital integrals.\p

We'd like to interpret their results in the case of very supercuspidal representations. The definition relies on the construction of supercuspidal representations by J.-K. Yu \cite{Yu01}. The construction begins with a datum consists of a sequence of subgroups $\G^0\subsetneq\G^1\subsetneq...\subsetneq\G^d=\G$ where each $\G^i$ is a Levi subgroup of $\G$ after base change to a tame extension of $F$, a depth-zero supercuspidal representation $\pi_0$ of $G^0=\G^0(F)$, and a sequence $(\phi^i)_{i=0}^d$ of characters on $(G^i)_{i=0}^d$. A supercuspidal representation is then constructed as the compact induction, from a compact open subgroup determined by $(\G^i)_{i=0}^d$ and depths of $(\phi_i)_{i=0}^d$, of some finite dimensional representation determined by $\pi_0$ and $(\phi_i)_{i=0}^d$. \p

Such supercuspidal representations were shown to exhaust all supercuspidal representations when $\text{char}(F)=0$ and $\text{char}(k)$ is large enough. If $\G^{d-1}$ is anisotropic (i.e. $G^{d-1}$ is compact, possibly modulo center), the supercuspidal representations constructed are called very supercuspidal representations. The result of Kim and Murnaghan, when restricted to a very supercuspidal representation $\pi$, says that there exists $T\in\LG$ satisfying $(C_{\G}(T))^o=\G^{d-1}$ such that locally near the identity, $\Theta_{\pi}(-)=\text{deg}(\pi)\hat{\mu}_T(-)$, where $\text{deg}(\pi)\in\mathbb{Q}$ is the formal degree of $\pi$.\p

On the other hand, the Lie algebra version of a theorem of Shalika \cite{Sh72} states that for any $T\in\LG$, there exists a lattice $\Lambda\subset\LG$ such that
\begin{equation}\label{preShalika}\mu_T(f)=\sum_{\CO\in\CO(0)}\Gamma_{\CO}(T)\mu_{\CO}(f),\;\forall f\in C_{c}^{\infty}(\LG/\Lambda).
\end{equation}

where $\CO(0)$ is the set of nilpotent orbits, and $\Gamma_{\CO}(X)$ are constants that depends on $T\in\LG$. These $\Gamma_{\CO}(X)$ are called the {\bf Shalika germs}. When $\Theta_{\pi}(-)=\text{deg}(\pi)\hat{\mu}_T(-)$, comparing (\ref{HCHlocal}) and (\ref{preShalika}) gives us  $c_{\CO}(\pi)=\text{deg}(\pi)\Gamma_{\CO}(T)$.\p

Our main goal in this article is to present an algorithm to compute these numbers $\Gamma_{\CO}(T)$. This gives the local character expansion at the identity for very supercuspidal representations as explained above. Our method begins with using test functions constructed via DeBacker's parametrization of nilpotent orbits \cite{De02b}. We first take $f$ to be the characteristic function of some nice set. The set has the property that it meets a nilpotent orbits $\CO$ and does not meet any other nilpotent orbits $\CO'$ of smaller or equal dimension. The test functions are then $f_{\pi^{2n}}(X)=f(\pi^{2n}X)$. Nilpotent orbital integrals satisfy the property $\mu_{\CO}(f_{\pi^{2n}})=q^{\dim\CO}\mu_{\CO}(f)$. This gives the term $\Gamma_{\CO}(T)\mu_{\CO}(f_{\pi^{2n}})$ a unique scaling behavior with respect to $n$. To compute $\Gamma_{\CO}(T)$, the essential thing remained is then to compute the asymptotic of $\mu_T(f_{\pi^{2n}})$.\p

To do so, we compute the orbital integral using the Cartan decomposition. We'll choose (based on $T$ and the test functions) specific points $x,y\in\mathcal{B}(\G,F)$ be two points in the Bruhat-Tits building of $\G$ and $G_x$, $G_y$ associated parahoric subgroups. Let $\mathbf{S}\subset\G$ be a maximal split torus whose corresponding apartment $\mathcal{A}=\mathcal{A}(\mathbf{S},F)$ contains $x$ and $y$. Denote by $\til{W}$ be the affine Weyl group of $\G$ (associated to $\mathbf{S}$) and $W_x, W_y\subset\til{W}$ the finite subgroups stabilizing $x$ and $y$, respectively. The Cartan decomposition is the following decomposition: $G=\bigsqcup_{w\in W_y\bsl\til{W}/W_x}G_ywG_x$.\p

Since we are doing orbital integral of $T$ who has a compact stabilizer, up to a normalizing constant we can write
$$\mu_T(f)=\int_Gf(\Ad(g)T)dg=\sum_{w\in W_y\bsl\til{W}/W_x}\int_{G_y}\int_{G_x}f(\Ad(g_y)\Ad(w)\Ad(g_x)T)dg_xdg_y.$$

For all but finitely many $w$ the integral will be zero (Theorem \ref{ellip}). We thus obtain $\mu_T(f)$ as a somewhat large combinatorial sum. We then prove some transversality result (Theorem \ref{transv}), in a spirit very close to that of the work of Goresky, Kottwitz and MacPherson \cite[Theorem 0.2 and Sec. 3.7]{GKM06}. The transversality result reduces each term in the sum into counting points on a quasi-finite cover of {\it Hessenberg varieties} (see \cite{GKM06} for definition of Hessenberg varieties) defined over the residue field $k$.\p

The combinatorics of the sum is a priori very complicated. Note the terms in the sum are indexed by $w\in W_y\bsl\til{W}/W_x$. Let $\CV=X_*(\mathbf{S})\otimes\mathbb{R}$ the real span of the cocharacters of $\mathbf{S}$. We have a map $ev:W_y\bsl\til{W}/W_x\ra\CV/W_y$ by $w\mapsto x-w^{-1}y$. Here $W_y$ acts linearly on $\CV$ (i.e. fixing the origin). For a certain real number $r$ related to our orbital integral problem, we'll draw walls (affine hyperplanes) of the form $\alpha(v)=0$ and $\alpha(v)=r$ on $\CV$, where $\alpha\in\Phi(\G,\mathbf{S})$ runs over roots of $\G$ with respect to $\mathbf{S}$. Since the set of these walls are invariant under $W_y$-action, they decompose $\CV/W_y$ into finitely many polyhedrons. Call this set of polyhedrons $\bar{\CP}$. We then have induced map $ev^{new}:W_y\bsl\til{W}/W_x\ra\bar{\CP}$.\p

These walls are where the behavior of the term indexed by $w$ changes. It turns out that we should group terms indexed by elements in $W_y\bsl\til{W}/W_x$ via their image in $\bar{\CP}$. In Theorem \ref{ellip} we prove that if $ev^{new}(w)$ is an unbounded polyhedron then the term indexed by such $w$ vanishes. Let $\bar{\CP}^{bd}\subset\bar{\CP}$ be the subset of bounded polyhedrons. We end up writing the orbital integral into

$$\mu_T(\text{test function})=\sum_{\Pi\in\bar{\CP}^{bd}}\sum_{w\in\Pi}(\text{the term indexed by }w).$$\p

Intuitively, what we then arrive is that there exists combinatorial constants $p_{\Pi}$ and geometric constant $J(\Pi)$ such that

$$\mu_T(\text{test function})=\sum_{\Pi\in\bar{\CP}^{bd}}p_{\Pi}\cdot J(\Pi).$$
and the corresponding result on Shalika germs
$$\Gamma_{\CO}(T)=\sum_{\Pi\in\bar{\CP}^{bd}}c_{p_{\Pi},\CO}\cdot J_{\CO}(\Pi).$$

See (\ref{main2}) and Theorem \ref{mainthm} for a precise statement. As mentioned, the numbers in the sum consist of the combinatorial part $c_{p_{\Pi},\CO}$ and the geometric part $J_{\CO}(\Pi)$. The latter part counts number of rational points on certain varieties over $k$ and is in general a mystery. Nevertheless, in Theorem \ref{degree} and Theorem \ref{vanish} we give conditions for $c_{p_{\Pi},\CO}$ and $J_{\CO}(\Pi)$ to vanish, respectively. This greatly reduces the complexity of the combinatorics. Roughly speaking, the complexity now is exponential in the difference in dimension of $\CO$ and that of a regular orbit. This makes it feasible to compute $\Gamma_{\CO}(T)$ for the top orbits for a given type of groups.\p

Our method gives for each $T$ and nilpotent orbit $\CO$ an explicit list of quasi-finite covers of Hessenberg varieties, for which the point-counting gives the geometric numbers $J_{\CO}(\Pi)$. These varieties are described in terms to the root system. In \cite[Sec. 3]{Ts15c} the author describes some examples of such varieties. For example, when $\CO$ comes from the top four nilpotent geometric orbits of a ramified unitary group and $T$ is certain nice half-integral depth elements, the author shows that $J_{\CO}(\Pi)$ and thus $\Gamma_{\CO}(T)$ appear to be numbers of points on varieties over $k$ whose $\ell$-adic cohomologies are generated by that of specific hyperelliptic curves.\p

We now briefly describe the structure of this article. In Section 2, we begin by introducing the notion of very elliptic elements. They are elements in the Lie algebra that arise from local characters of very supercuspidal representations. After proving such elements have anisotropic stabilizer (Theorem \ref{ellip}) and the transversality result (Theorem \ref{transv}) in 2.1, we begin with a general discussion of orbital integrals of a very elliptic element $\til{\phi}_x$. In 2.2 we apply the Cartan decomposition and eventually in 2.3 write the orbital integral into a sum of products of combinatorial number and geometric numbers in (\ref{main}).\p

In Section 3 we turn to the discussion of Shalika germs. In 3.1 we describe how to find a sequence of suitable test functions using DeBacker's parametrization of nilpotent orbit, and explain how we use the homogeneity property of nilpotent orbital integrals. We then in 3.2 describe how the notions in 2.3 can be used to derive our main result for Shalika germ, which is Theorem \ref{mainthm}. Lastly, in 3.3, we prove a vanishing result that greatly simplifies the formula in practice.\p

\subsection*{Acknowledgments}The author would like to thank his Ph.D. advisor Benedict Gross for all inspiration and encouragement. He would like to thank Zhiwei Yun for many beneficial discussions, from which in particular the method in the first half of this article was improved. He would also like to thank Thomas Hales, Tasho Kaletha, Bao Le Hung, Loren Spice, Jack Thorne, Pei-Yu Tsai and Jerry Wang for helpful suggestions and discussions.\p

\section{A computation for orbital integrals}

\subsection{Very elliptic elements}\label{orb}

Our purpose is to compute some orbital integral on the $\LG$. Let $\G^{sc}$ be the simply connected cover of the derived group $\G^{der}$ of $\G$. $\G^{sc}$ is also the simply connected cover of the adjoint quotient $\G^{ad}$ of $\G$. Since we assume $\text{char}(F)$ is very good for $\G$, in the exact sequence $\G^{sc}(F)\ra\G^{ad}(F)\ra H^1(F,Z(\G^{sc}))$ the last term is finite. \p

The map $\G^{sc}\ra\G^{ad}$ factors through $\G^{sc}\ra\G\ra\G^{ad}$. Therefore the image of $\G^{sc}(F)$ in $\G^{ad}(F)$ is a finite index subgroup of the image of $G=\G(F)$. When $\text{char}(F)$ is very good for $\G$, the maps $\G^{sc}\ra\G\ra\G^{ad}$ induce full-rank differentials on the Lie algebras. Computation of orbital integrals on the Lie algebra can thus be reduced from $\G$ to $\G^{sc}$. We thus assume from now on, as mentioned in the introduction, that $\G=\G^{sc}$ is semi-simple and simply connected.\p

In this and the following section, beside the group $\G$ over the local field $F$, $x,y$ will be two points on the building $\mathcal{B}(\G,F)$ of $\G$ over $F$. We will assume that $\G$ is tamely ramified, i.e. $\G$ splits over a tamely ramified extension. We also assume that the coordinate of $x$ has denominator coprime to $p$, i.e. $x$ becomes a hyperspecial point after base change to a tame extension. Take $d_x>d_y$ two real numbers. For our computation we'll make use of a maximal split torus $\mathbf{S}$ whose corresponding apartment $\mathcal{A}=\mathcal{A}(\mathbf{S},F)$ contains $x$ and $y$. We fix such an $\mathbf{S}$ from now on. Let $G_x$ and $G_y$ be the parahorics stabilizing $x$ and $y$ respectively in $G$.\p

Let $\til{\mathbb{R}}$ be the totally ordered set of symbols $\til{\mathbb{R}}=\{r,r+\,|\,r\in\mathbb{R}\}$ and likewise for $\til{\mathbb{R}}_{\ge 0}$. We have Moy-Prasad filtrations \cite[Sec. 2]{MP94} $\{G_{x,r}\}_{r\in\til{\mathbb{R}}_{\ge 0}}$ and $\{\LG_{x,r}\}_{r\in\til{\mathbb{R}}}$, where $G_{x,r+}:=\lim_{s\ra r+}G_{x,s}$ and $\LG_{x,r+}:=\lim_{s\ra r+}\LG_{x,s}$.
We'll denote by $\LG_{x,r:s}=\LG_{x,r}/\LG_{x,s}$ for any $r<s$ in $\til{\mathbb{R}}$. Finally define $L_x=G_x/G_{x,0+}$ and $V_x=\LG_{x,d_x:d_x+}$. The same notations will be used when $x$ is replaced by $y$ or other points on the building.\p

Since we assume $\G$ to be simply connected, the group $L_x$ is the $k$-points of a connected reductive group $\mathbf{L}_x$ defined over $k$. $L_x$ acts on the $k$-vector space $V_x$ by adjoint action $\text{Ad}$, and this action also arises from an algebraic representation of $\mathbf{L}_x$. We consider an element $\phi_x\in V_x$ with the following assumptions: (1) the $\mathbf{L}_x$ orbit of $\phi_x$ is (Zariski) closed and (2) The stabilizer of $\phi_x$ in $\mathbf{L}_x$ is an anisotropic torus.\p

Let $\til{\phi}_x\in\LG_{x,d_x}$ be any lift\footnote{This is what we meant by a very elliptic element in the title of this subsection.} of $\phi_x$, and $\phi_y\in V_y$ be any element. The main goal of this section is to compute the following orbital integral:

\begin{equation}\label{orbit}\mu_{\til{\phi}_x}(1_{\phi_y+\LG_{y,d_y+}})=\;?
\end{equation}

Here $1_{\phi_y+\LG_{y,d_y+}}$ denote the characteristic function of $\phi_y+\LG_{y,d_y+}$ and $\mu_{\til{\phi}_x}$ the orbital integral of $\til{\phi}_x$.
\p

We shall explain more about this setting, especially about the point $x$, the group $\mathbf{L}_x$ and its representation $V_x$. For the result in this subsection we have to assume $\G$ is tamely ramified. We first recall some results in \cite[Sec. 4]{RY14}. Let's say $x$ has denominator $m$. Then $md_x$ must be an integer, since $V_x=0$ otherwise. Fix $\zeta_m$ a $m$-th root of unity in $\bar{k}$. By \cite[Theorem 4.1]{RY14}, there is an algebraic group $\mathsf{G}$ defined over $k$ together with an order $m$ automorphism $\theta:\mathsf{G}\ra\mathsf{G}$, which is only defined over $k(\zeta_m)$, but such that the grading $\text{Lie }\mathsf{G}=\bigoplus_{i\in\Z/m}(\text{Lie }\mathsf{G})(i)$, where $(\text{Lie }\mathsf{G})(i):=\{X\in\text{Lie }\mathsf{G}\,|\,\theta(X)=\zeta_m^{i}\}$, is defined over $k$. In particular the fixed subgroup $\mathsf{G}^{\theta}$ is also a $k$-subgroup of $\mathsf{G}$.\p

The statement is that the automorphism $\theta$ is such that we have compatible $k$-isomorphisms $\mathbf{L}_x\cong(\mathsf{G}^{\theta})^o$ and $V_x\cong(\text{Lie }\mathsf{G})(md_x)$, that is, the action of $\mathbf{L}_x$ on $V_x$ agrees with the adjoint action of $(\mathsf{G}^{\theta})^o$ on $(\text{Lie }\mathsf{G})(md_x)$. We can therefore think of $\phi_x$ as an element in $(\text{Lie }\mathsf{G})(md_x)$ with closed orbit and anisotropic stabilizer under the action of $(\mathsf{G}^{\theta})^o$. In addition, that $\phi_x$ has closed orbit is equivalent to $\phi_x\in\text{Lie }\mathsf{G}$ is semisimple \cite[Lemma 2.12 and Cor. 2.13]{Le09}. In this subsection we prove two crucial results for the orbital integral of $\til{\phi}_x$, that is Theorem \ref{ellip} and Theorem \ref{transv}:\p

\begin{theorem}\label{ellip} $\til{\phi}_x$ has anisotropic stabilizer in $\G$.
\end{theorem}
\begin{proof}For any $g\in G$ which centralizes $\til{\phi}_x$, we shall prove that $g\in G_x$. Supoose otherwise that $g.x=x'\not=x$. Then since $g$ centralizes $\til{\phi}_x$ we have $\til{\phi}_x\in\LG_{x',d_x}$ as well.\p

Without loss of generality we may assume $x'\in\mathcal{A}=\mathcal{A}(\mathbf{S},F)$. Let $\Psi(\G,\mathbf{S})$ be the set of affine roots of $\mathbf{S}$ in $\G$. One has
$$V_x=\bigoplus_{\psi\in\Psi(\G,\mathbf{S}),\,\psi(x)=d_x}\LG_{x,d_x:d_x+,\dot{\psi}}\;,\text{ and}$$
$$\phi_x\in\text{Im}(\LG_{x,d_x}\cap\LG_{x',d_x}\ra V_x)=\bigoplus_{\psi\in\Psi(\G,\mathbf{S}),\,\psi(x')\ge\psi(x)=d_x}\LG_{x,d_x:d_x+,\dot{\psi}}.$$\pp

where $\dot{\psi}$ denotes the gradient of $\psi$, and $\LG_{x,d_x:d_x+,\dot{\psi}}$ is the graded piece in the root space $\LG_{\dot{\psi}}$ of depth $d_x$. Since we assume $x$ has rational coordinates with denominator $m$, $m(x'-x)\in X_*(\mathbf{S})$ corresponds to a cocharacter $\lambda:\mathbb{G}_m\ra\mathbf{L}_x$ with image in the maximal $k$-split torus of $\mathbf{L}_x$ that corresponds to $\mathbf{S}$. The previous discussion says $\phi_x\in V_x$ lies in the non-negative weight space of $\lambda$.\p

Identify $\mathbf{L}_x$ as $(\mathsf{G}^{\theta})^o$ and $V_x$ with $(\text{Lie }\mathsf{G})(md_x)$. Since $\lambda$ has image in $\mathsf{G}^{\theta}$, the $\Z$-grading given by $\lambda$ is compatible with the grading given by $\theta$, i.e. we can write
$$\text{Lie }\mathsf{G}=\bigoplus_{i\in\Z/m,\,j\in\Z}(\text{Lie }\mathsf{G})(i)_j$$
where $(\text{Lie }\mathsf{G})(i)_j$ is the subspace of $(\text{Lie }\mathsf{G})(i)$ on which $\lambda$ acts by $z^j$. We have $\phi_x\in(\text{Lie }\mathsf{G})(md_x)_{\ge 0}$ is semisimple. $\mathsf{G}$ has a parabolic subgroup $\mathsf{P}$ with Lie algebra $\Lie\mathsf{P}=(\Lie\mathsf{G})_{\ge 0}$ of non-negative weights of $\lambda$, and $\mathsf{P}$ has a Levi subgroup $\mathsf{L}$ with $\Lie\mathsf{L}=(\Lie\mathsf{G})_0$. Denote by $\phi_{x,0}$ the image of $\phi_x$ to $\Lie\mathsf{L}$ which is just given by projection to the weight-$0$ subspace. \p

By the action of $\lambda(z)$ with $z\ra 0$, we see $\phi_{x,0}$ is in the (Zariski) closure of the $(\mathsf{G}^{\theta})^o$-orbit of $\phi_x$ and thus $\phi_x$ is conjugate to $\phi_{x,0}$ by elements in $(\mathsf{G}^{\theta})^o(\bar{k})$. One can possibly compute Galois cohomology (especially as groups are now defined over the finite field $k$) to show that $\phi_x$ is conjugate to $\phi_{x,0}$ by elements in $(\mathsf{G}^{\theta})^o(k)$. If so, then $\phi_{x,0}$ should have anisotropic stabilizer in $(\mathsf{G}^{\theta})^o(k)$ but $\lambda(\mathbb{G}_m)$ centralizes $\phi_{x,0}$, and we derive the contradiction we want. However we are not able to carry out a general Galois cohomology argument. Instead, we prove the following lemma:\p

\begin{lemma}\label{unip} Let $\mathsf{U}$ be the unipotent radical of $\mathsf{P}$ with $\Lie\mathsf{U}=(\Lie\mathsf{G})_{>0}$. Then for $k$-rational $X\in(\Lie\mathsf{G})_0$ and $Y\in(\Lie\mathsf{G})_{>0}$ such that $X$ and $X+Y$ are semisimple, there exists $u\in\mathsf{U}(k)$ satisfying $\Ad(u)(X+Y)=X$.\pp

Moreover, if $X,Y\in(\Lie\mathsf{G})(md_x)$, we may take $u\in(\mathsf{G}^{\theta})^o$.
\end{lemma}\pp
\begin{proof} We may assume $Y\in(\Lie\mathsf{G})_{\ge n}\bsl(\Lie\mathsf{G})_{>n}$ for some integer $n>0$. It suffices to find $u\in\mathsf{U}(k)$ such that $\Ad(u)(X+Y)\equiv X$ modulo $(\Lie\mathsf{G})_{>n}$ then use induction on $m$, as $(\Lie\mathsf{G})_{>n}=0$ for $m\gg0$. Write $Y=Y_n+Y_{>n}$ with $Y_m\in(\Lie\mathsf{G})_{n}(k)$ and $Y_{>n}\in(\Lie\mathsf{G})_{>n}(k)$ (note $(\Lie\sG)_n(k)$ means $k$-points of $\Lie\sG$, same for all other places of this proof). If we can find $V\in(\Lie\mathsf{G})_{n}(k)$ such that $[V, X]=-Y_n$, then $\Ad(\se(V))(X+Y)-X\in(\Lie\mathsf{G})_{>n}$ and we are done. Note here if $\text{char}(k)$ is not large enough we cannot have $\se$ as the naive exponential map $\se$; we'll not deal with this issue until the end. Now it suffices to prove that the image of $\text{ad}(X):(\Lie\mathsf{G})_n\ra(\Lie\mathsf{G})_n$ contains $Y_n$. This statement is independent of the base field and we work over an algebraic closure $\bar{k}$.\p

Take a maximal torus $\mathsf{T}$ of $\mathsf{G}$ that contains $X$ and is contained in $C_{\mathsf{G}}(\lambda(\mathbb{G}_m))$. We claim that if $X+Y$ is semisimple, then for every root $\alpha\in\Phi(\mathsf{G},\mathsf{T})$ with $\langle\lambda,\alpha\rangle=n$, we have $(d\alpha)(X)\not=0$ or $Y_{n,\alpha}=0$. Assuming the claim, we have
$$ad(X)(\sum_{\alpha}[(d\alpha)(X)]^{-1}Y_{n,\alpha})=Y_n.$$

The claim is a statement only on a particular root space, and can be checked directly by embedding $\mathsf{G}$ into $\text{GL}_N$ and $\mathsf{T}$ into the diagonal torus. This proves the lemma except for the last sentence. For the last sentence, note that when we take $V\in(\Lie\mathsf{G})_{n}$ such that $[V, X]=-Y_n$, since $Y_n\in(\Lie\mathsf{G})(md_x)_n$, it makes no harm to project $V$ to its component in $(\Lie\mathsf{G})(0)_{n}$. One then has $\se(V)\in(\mathsf{G}^{\theta})^o$.\p

Lastly, as promised we justify the exponential map $\se$ on $(\Lie\sG)_n(k)\subset(\Lie\sG)_{>0}(k)$ when $\text{char}(k)$ is small by finding an approximation. Let $\Phi^+=\{\alpha\in\Phi(\sG,\mathsf{T}),\,\langle\alpha,\lambda\rangle>0\}$ be the subset with positive weight with respect to $\lambda$, so $\Lie\sU=\bigoplus_{\alpha\in\Phi^+}(\Lie\sU)_{\alpha}$ (this decomposition might not be defined over $k$). On each $(\Lie\sU)_{\alpha}$, $\alpha\in\Phi^+$ we have the exponential map $\se_{\alpha}:(\Lie\sU)_{\alpha}\cong\sU_{\alpha}$, which is an isomorphism to the root subgroup $\sU_{\alpha}$. Taking product in any order gives an isomorphism $\se_n:(\Lie\sU)_n\cong\sU_n:=\prod_{\langle\alpha,\lambda\rangle\ge n}\sU_{\alpha}/\prod_{\langle\alpha,\lambda\rangle>n}\sU_{\alpha}$ defined over $k$ (which does not depend on the order for the product).\p

Now we can define a map on $k$-points of $\Lie\sU$ to $\sU$ compatible with all $\se_n$ following \cite[Sec. 1.3]{Ad98}. This is done by writing an element in $V\in(\Lie\sU)(k)$ as $V=\sum_{i\ge 1}V_i$, $V_i\in(\Lie\sU)_i(k)$ and take $\se(V)=\se(V_1)\cdot\se(V_2)\cdot...$ where $\se(V_i)$ are chosen representatives in $\prod_{\langle\alpha,\lambda\rangle\ge i}\sU_{\alpha}$ that are defined over $k$. We can also choose $\theta$-fixed representatives whenever $V_i\in(\Lie\sG)(0)$, which is needed for the last sentence of the lemma.
\end{proof}\p

Lemma \ref{unip} gives explicit rational element $u\in(\mathsf{G}^{\theta})^o(k)$ such that $\Ad(u)(\phi_x)=\phi_{x,0}$, and thus $\phi_{x,0}$ has anisotropic stabilizer and we have a contradiction unless $\lambda$ is trivial, that is $x'=x$. This finishes the proof of Theorem \ref{ellip}.\end{proof}\pp

With Theorem \ref{ellip}, up to a normalizing constant we have
$$\mu_{\til{\phi}_x}(1_{\phi_y+\LG_{y,d_y+}})=\int_{G_x}1_{\phi_y+\LG_{y,d_y+}}(\Ad(g)\til{\phi}_x)dg.$$\p

\begin{theorem}\label{transv} The value of $\int_{G_x}1_{\phi_y+\LG_{y,d_y+}}(\Ad(g)\til{\phi}_x)dg$ depends only on $\phi_x$ (but not on the lift $\til{\phi}_x$ chosen).
\end{theorem}
\begin{proof}We may decompose $G_x=\bigsqcup_{\dot{g}G_{x,0+}}\dot{g}G_{x,0+}$ into left $G_{x,0+}$-cosets, and it suffices to prove that for every $\dot{g}$

\begin{equation}\label{transv1}\int_{G_{x,0+}}1_{\Ad(\dot{g}^{-1})(\phi_y+\LG_{y,d_y+})}(\Ad(h)\til{\phi}_x)dh\;\text{ depends only on }\phi_x.\end{equation}\pp\pp

By letting $y=\dot{g}^{-1}y$ and $\phi_y=\Ad(\dot{g}^{-1})(\phi_y)$ we may drop $\dot{g}$ in (\ref{transv1}). We claim the following general statement: Let $r\in\til{\mathbb{R}}_{>0}$ and $\til{\phi}_x\in\phi_x+\LG_{x,d_x+}$ be any element in the coset. Then
\begin{equation}\label{transv2}\int_{G_{x,r}}1_{\phi_y+\LG_{y,d_y+}}(\Ad(h)\til{\phi}_x)=\frac{|G_{x,r}|\cdot|(\til{\phi}_x+\LG_{x,d_x+r})\cap(\phi_y+\LG_{y,d_y+})|}{|(\til{\phi}_x+\LG_{x,d_x+r})|}.
\end{equation}\p

Note the measures $|(\til{\phi}_x+\LG_{x,d_x+r})|=|\LG_{x,d_x+r}|$ and $|(\til{\phi}_x+\LG_{x,d_x+r})\cap(\phi_y+\LG_{y,d_y+})|=|\LG_{x,d_x+r}\cap\LG_{y,d_y+}|$ if the first intersection is non-empty and zero otherwise. Taking $r=0+$ in (\ref{transv2}) gives (\ref{transv1}). To prove (\ref{transv2}) we use induction on $r$ in the way that we assume it works for all $r'>r$ and prove the statement for $r$. This works because $G_{x,r}$ jumps semicontinuously, and that when $r$ is large enough, $\LG_{x,d_x+r}\subset\LG_{y,d_y+}$ in which case the conclusion is immediate.\p

(\ref{transv2}) is clear when $r$ is not a jump for the Moy-Prasad filtration of $G_x$ or when $(\til{\phi}_x+\LG_{x,d_x+r})\cap(\phi_y+\LG_{y,d_y+})=\emptyset$. Suppose otherwise $(\til{\phi}_x+\LG_{x,d_x+r})\cap(\phi_y+\LG_{y,d_y+})\not=\emptyset$ and denote by $r'$ is the next jump of the filtration, by comparing integration over $G_{x,r}$ and over $G_{x,r'}$, the induction is exactly to show that
\begin{equation}\label{cosetcount}\#\{\dot{h}\in G_{x,r}/G_{x,r'}\,|\,(\Ad(\dot{h})\til{\phi}_x+\LG_{x,d_x+r'})\cap(\phi_y+\LG_{y,d_y+})\not=\emptyset\}
\end{equation}
$$=\frac{|G_{x,r}|\cdot|\LG_{x,d_x+r'}|\cdot|\LG_{x,d_x+r}\cap\LG_{y,d_y+}|}{|G_{x,r'}|\cdot|\LG_{x,d_x+r}|\cdot|\LG_{x,d_x+r'}\cap\LG_{y,d_y+}|}.$$\p

We may decompose
$$\LG_{x,d_x+r:d_x+r'}=\bigoplus_{\alpha(x-y)<d_x+r-d_y}\LG_{x,d_x+r:d_x+r',\alpha}\oplus\bigoplus_{\alpha(x-y)\ge d_x+r-d_y}\LG_{x,d_x+r:d_x+r',\alpha}.$$\p

Let's now in the remaining proof of this lemma denote the first summand by $W_1$ and the second by $W_2$. They are subspaces of $\LG_{x,d_x+r:d_x+r'}$ as $k$-vector spaces. In the decomposition, $W_1$ is the image of $\LG_{x,d_x+r}\cap\LG_{y,d_y+}$, consequently the RHS of (\ref{cosetcount}) is equal to $\frac{\#(G_{x,r}/G_{x,r'})}{\#W_2}.$\p

Given the assumption $(\til{\phi}_x+\LG_{x,d_x+r})\cap(\phi_y+\LG_{y,d_y+})\not=\emptyset$, that whether $(\Ad(\dot{h})\til{\phi}_x+\LG_{x,d_x+r'})\cap(\phi_y+\LG_{y,d_y+})\not=\emptyset$ depends only on the image of $\Ad(\dot{h})\til{\phi}_x-\til{\phi}_x$ in $\LG_{x,d_x+r:d_x+r'}$, and actually only on its projection to $W_2$. In fact, let $\pi_{W_2}$ be the projection to from $\LG_{x,d_x+r:d_x+r'}$ or $\LG_{x,d_x+r}$ to $W_2$. Then
$$(\Ad(\dot{h})\til{\phi}_x+\LG_{x,d_x+r'})\cap(\phi_y+\LG_{y,d_y+})\not=\emptyset
\text{ and }
(\Ad(\dot{h}')\til{\phi}_x+\LG_{x,d_x+r'})\cap(\phi_y+\LG_{y,d_y+})\not=\emptyset$$
$$\Rightarrow\Ad(\dot{h})\til{\phi}_x-\Ad(\dot{h}')\til{\phi}_x\in\LG_{x,d_x+r'}+\LG_{y,d_y+}
\Rightarrow\pi_{W_2}(\Ad(\dot{h})\til{\phi}_x-\Ad(\dot{h}')\til{\phi}_x)=0,$$

that is, the intersection is non-empty exactly when $\pi_{W_2}(\Ad(\dot{h})\til{\phi}_x-\til{\phi}_x)$ hits a specific element in $W_2$.
\p

We make use of the mock-exponential map $\se(\cdot):\LG_{x,0+}\ra G_{x,0+}$ defined in \cite[Appendix A]{AS09}, which is a bijection satisfying properties as if the exponential map was defined. In particular we identify $\se:\LG_{x,r}/\LG_{x,r'}\xra{\sim}G_{x,r}/G_{x,r'}$. Observe that $\Ad(\dot{h})\til{\phi}_x+\LG_{x,d_x+r'}=\til{\phi}_x+[\se^{-1}(\dot{h}),\phi_x]+\LG_{x,d_x+r'}$.\p

\begin{claim}\label{ellipsurj}The linear transformation (over $k$) from $\LG_{x,r:r'}$ to $W_2$ given by $X\mapsto\pi_{W_2}([X,\phi_x])$ is surjective.
\end{claim}

Suppose the validity of the claim. Then $\pi_{W_2}(Ad(\dot{h})\til{\phi}_x-\til{\phi}_x)$ runs over all elements in $W_2$. Thus  $(\Ad(\dot{h})\til{\phi}_x+\LG_{x,d_x+r'})\cap(\phi_y+\LG_{y,d_y+})\not=\emptyset$ exactly $(\#W_2)^{-1}$ of the time. This is exactly (\ref{cosetcount}).\p

It remains to prove Claim \ref{ellipsurj}. By \cite[Proposition 4.1]{AR00}, as we assume $\text{char}(k)$ is very good, there is a $G$-invariant non-degenerate symmetric bilinear form $\kappa(\cdot,\cdot)$ on $\LG$ whose reduction $\bar{\kappa}$ identifies $\LG_{x,d_x+r:d_x+r'}$ with the dual of $\LG_{x,-d_x-r:(-d_x-r)+}$. For every root $\alpha\in\Phi(\G,\mathbf{S})$, it identifies $\LG_{x,d_x+r:d_x+r',\alpha}$ with $\LG_{x,-d_x-r:(-d_x-r)+,-\alpha}$. In particular, we can identify the dual of $W_2$:

$$W_2^*=\bigoplus_{\alpha(x-y)\le d_x+r-d_y}\LG_{x,-(d_x+r):(-d_x-r)+,\alpha}\subset\LG_{x,-d_x-r:(-d_x-r)+}$$\pp\pp

Suppose the map in Claim \ref{ellipsurj} is not surjective. Then there exists a non-zero element $N\in W_2^*\subset\LG_{x,-d_x-r:(-d_x-r)+}$ such that $\bar{\kappa}([X,\phi_x],N)=0\in k$ for all $X\in\LG_{x,r:r+}$. However $\bar{\kappa}([X,\phi_x],N)=\bar{\kappa}(X,[\phi_x,N])$. Hence we have $[\phi_x,N]=0$.\p

Recall we have algebraic group $\sG$ with an order $m$ automorphism $\theta$ which gives $\Z/m$-grading $\Lie\sG=\bigoplus_{i\in\Z/m}(\Lie\sG)(i)$ such that we may realize $\phi_x\in(\Lie\sG)(md_x),\; N\in(\Lie\sG)(m(-d_x-r))$ with $[\phi_x,N]=0$. Note $N$ is a nilpotent element, and we shall derive a contradiction with the assumption that $\phi_x$ has anisotropic stabilizer in $(\sG^{\theta})^o$.\p

Let $\mathsf{H}:=C_{\sG}(\phi_x)^o$ be the connected centralizer of $\phi_x$. By assumption $(\mathsf{H}^{\theta})^o$ is an anisotropic torus. Since $\theta$ acts on $\phi_x$ by a constant, $\mathsf{H}$ is invariant under $\theta$ and we have likewise a $\Z/m$-grading on $\Lie\mathsf{H}$. In particular $N\in(\Lie\mathsf{H})(m(-d_x-r))$. The proof in \cite[Lemma 2.11]{Le09} gives a unique one-parameter subgroup in $(\mathsf{H}^{\theta})^o$ which acts non-trivially on the line generated by $N$. But this contradicts with the assumption $(\mathsf{H}^{\theta})^o$ is anisotropic.\end{proof}\p

\subsection{Cartan decomposition}

We are now ready to evaluate (\ref{orbit}). Because of Theorem \ref{ellip} and that $\phi_y+\LG_{y,d_y+}$ is invariant under $G_{y,0+}$, up to a normalizing constant, one can write
$$\mu_{\til{\phi}_x}(1_{\phi_y+\LG_{y,d_y+}})=\int_{G_{y,0+}\bsl G}1_{\phi_y+\LG_{y,d_y+}}(\Ad(g)\til{\phi}_x)dg.$$

We have the Cartan decomposition $G_y\bsl G/G_x\cong W_y\bsl\til{W}/W_x$, where $\til{W}$ is as in \cite[1.2]{Tits} and $W_x$ and $W_y$ are the stabilizers of $x$ and $y$ in $\til{W}$, respectively. We shall choose an arbitrary set of representatives for $W_y\bsl\til{W}/W_x$ in $G$; whenever we write $w\in W_y\bsl\til{W}/W_x$, $w$ can be understood as an element in $\til{W}$ or $G$ (while whatever follows should be essentially independent of this choice). Note that since we assume $\G$ to be semi-simple and simply connected, $\til{W}$ agrees with the affine Weyl group.\p

Suppose we normalize our measure so that $G_{y,0+}$ has measure $1$. With the Cartan decomposition we can rewrite the integral:

$$\int_{G_{y,0+}\bsl G}1_{\phi_y+\LG_{y,d_y+}}(\Ad(g)\til{\phi}_x)dg=\sum_{w\in W_y\bsl\til{W}/W_x}\sum_{g\in G_{y,0+}\bsl G_ywG_x}1_{\phi_y+\LG_{y,d_y+}}(\Ad(g)\til{\phi}_x)\s\s$$
$$=\sum_{w\in W_y\bsl\til{W}/W_x}\sum_{h\in L_y}\sum_{\dot{g}\in G_y\bsl G_ywG_x}1_{\phi_y+\LG_{y,d_y+}}(\Ad(h)\Ad(\dot{g})\til{\phi}_x)\;\;\;\,$$
$$\s\s=\sum_{w\in W_y\bsl\til{W}/W_x}\sum_{h\in L_y}\sum_{\dot{\gamma}\in(G_x\cap G_{w^{-1}y})\bsl G_x}1_{\Ad(w^{-1}h^{-1})(\phi_y)+\LG_{w^{-1}y,d_y+}}(\Ad(\dot{\gamma})\til{\phi}_x).$$\pp

Here the point $w^{-1}y\in\mathcal{A}$ on the building is characterized by $G_{w^{-1}y}=\Ad(w^{-1})G_y$, and $\Ad(w^{-1})$ sends $\Ad(h^{-1})(\phi_y)\in\LG_{y,d_y:d_y+}$ to $\Ad(w^{-1}h^{-1})(\phi_y)\in\LG_{w^{-1}y,d_y:d_y+}$.\p

We'll show that there are only finitely many $w\in W_y\bsl\til{W}/W_x$ such that $\Ad(g)\til{\phi}_x$ is in the compact lattice $\LG_{y,d_y}$, i.e. only finitely many $w$ in the sum contributes (see Lemma \ref{finite}). The order of sum and integral thus don't matter. For convenience write $\phi_{w^{-1}y}=\Ad(w^{-1})\phi_y\in\LG_{w^{-1}y,d_y:d_y+}$. By Theorem \ref{transv}, we can, by assigning $\LG_{x,d_x+}$ to have Haar measure $1$ in the following integral, rewrite 

$$\mu_{\til{\phi}_x}(1_{\phi_y+\LG_{y,d_y+}})\s\s\s\s\s\s\s\s\s\s\s\s$$
$$=\sum_{w\in W_y\bsl\til{W}/W_x}\sum_{h\in L_{w^{-1}y}}\sum_{\dot{\gamma}\in(G_x\cap G_{w^{-1}y})\bsl G_x}\int_{\phi_x+\LG_{x,d_x+}}1_{\Ad(h^{-1})(\phi_{w^{-1}y})+\LG_{w^{-1}y,d_y+}}(\Ad(\dot{\gamma})\til{\phi}_x)d\til{\phi}_x\s$$
$$=\sum_{w\in W_y\bsl\til{W}/W_x}q^{-m_w'}\sum_{h\in L_{w^{-1}y}}\sum_{\dot{\gamma}\in(G_x\cap G_{w^{-1}y})\bsl G_x}1_{(\Ad(h^{-1})(\phi_{w^{-1}y})+\LG_{w^{-1}y,d_y+}+\LG_{x,d_x+})/\LG_{x,d_x+}}(\Ad(\dot{\gamma})\phi_x).$$\pp\pp

where $m_w'$ is the non-negative number such that $q^{m_w'}=[\LG_{x,d_x+}:\LG_{x,d_x+}\cap\LG_{w^{-1}y:d_y+}]$. Now one note that the sum over $(G_x\cap G_{w^{-1}y})\bsl G_x$ is $[(G_x\cap G_{w^{-1}y})G_{x,0+}:G_x\cap G_{w^{-1}y}]$ times the sum over $(G_x\cap G_{w^{-1}y})G_{x,0+}\bsl G_x$, which is a quotient of $L_x$ by a parabolic subgroup.\p

More precisely, one knows \cite[3.5]{Tits} $(G_x\cap G_{w^{-1}y})G_{x,0+}\bsl G_x=\mathbf{P}_{x,w}(k)\bsl\mathbf{L}_x(k)$, for some algebraic subgroup $\mathbf{P}_{x,w}\subset\mathbf{L}_x$ as given in the following lemma. Note that by Lang's theorem $\mathbf{P}_{x,w}(k)\bsl\mathbf{L}_x(k)=(\mathbf{P}_{x,w}\bsl\mathbf{L}_x)(k)$.\pp\pp

\begin{lemma}\label{tech} Let $\bar{\mathbf{S}}$ be the maximal $k$-split torus of $\mathbf{L}_x$ that correspond to $\mathbf{S}$ as in \cite[3.5]{Tits}. In particular, the character lattice $X^*(\mathbf{S})$ can be identified to $X^*(\bar{\mathbf{S}})$. The group $\mathbf{P}_{x,w}$ is then the parabolic subgroup of $\mathbf{L}_x$ generated by $Z_{\mathbf{L}_x}(\bar{\mathbf{S}})$ and the root subgroups of those roots $\alpha\in\Phi(\mathbf{L}_x,\bar{\mathbf{S}})$ satisfying $\alpha(x-w^{-1}y)\le 0$, where $x-w^{-1}y\in X_*(\mathbf{S})\otimes\mathbb{R}$ is evaluated by $\alpha$ via the identification $X^*(\mathbf{S})\cong X^*(\bar{\mathbf{S}})$.
\end{lemma}\pp

We now rewrite the integral:
$$\mu_{\til{\phi}_x}(1_{\phi_y+\LG_{y,d_y+}})\s\s\s\s\s\s\s\s\s\s\s\s$$
$$=\sum_{w\in W_y\bsl\til{W}/W_x}q^{m_w-m_w'}\sum_{h\in L_{w^{-1}y}}\sum_{\dot{\gamma}\in(\mathbf{P}_{x,w}\bsl\mathbf{L}_x)(k)}1_{(\Ad(h^{-1})(\phi_{w^{-1}y})+\LG_{w^{-1}y,d_y+}+\LG_{x,d_x+})/\LG_{x,d_x+}}(\Ad(\dot{\gamma})\phi_x),$$\p

where $m_w$ is a non-negative integer such that $q^{m_w}=[(G_x\cap G_{w^{-1}y})G_{x,0+}:G_x\cap G_{w^{-1}y}]=[G_{x,0+}:G_{x,0+}\cap G_{w^{-1}y}]=[\LG_{x,0+}:\LG_{x,0+}\cap\LG_{w^{-1}y}]$.\p

\subsection{Organization of double cosets $W_ywW_x$} We shall rewrite the last equation in the previous subsection. Note that $\Ad(\dot{\gamma})(\phi_x)\in\LG_{x,d_x:d_x+}$. Thus we may replace $(\Ad(h^{-1})(\phi_{w^{-1}y})+\LG_{w^{-1}y,d_y+}+\LG_{x,d_x+})/\LG_{x,d_x+}$ by $((\Ad(h^{-1})(\phi_{w^{-1}y})+\LG_{w^{-1}y,d_y+}+\LG_{x,d_x+})\cap\LG_{x,d_x})/\LG_{x,d_x+}$. The latter is contained in (\ref{Hx}) of the following: 

\begin{equation}\label{Hx}((\LG_{w^{-1}y,d_y}+\LG_{x,d_x+})\cap\LG_{x,d_x})/\LG_{x,d_x+}
=\dsp\sum_{\alpha\in\Phi(\G,\mathbf{S}),\,\alpha(x-w^{-1}y)\le d_x-d_y}\LG_{x,d_x:d_x+,\alpha},
\end{equation}

where for any $r,s\in\til{\mathbb{R}}$ we define $\LG_{x,r,\alpha}$ to be the intersection $\LG_{x,r}\cap\mathfrak{u}_{\alpha}$ ($\mathfrak{u}_{\alpha}$ being the root space of $\alpha$), and $\LG_{x,r:s,\alpha}=\LG_{x,r,\alpha}/\LG_{x,s,\alpha}\subset\LG_{x,r:s}$.\p

Consider the vector space $\V:=X_*(\mathbf{S})\otimes\mathbb{R}$. We draw affine hyperplanes $\{v\;|\;\alpha(v)=0\}$ and $\{v\;|\;\alpha(v)=d_x-d_y\}$. These hyperplane cut the whole space into polyhedrons of various dimensions. We have\pp\pp

\begin{lemma}\label{finite}If $w\in W_y\bsl\til{W}$ is such that $x-w^{-1}y$ lies in an unbounded polyhedron, then $((\LG_{w^{-1}y,d_y}+\LG_{x,d_x+})\cap\LG_{x,d_x})/\LG_{x,d_x+}
\subset\LG_{x,d_x:d_x+}$ contains no element with closed $\mathbf{L}_x$-orbit and anisotropic stabilizer.
\end{lemma}

This implies that such $w$, or say their image in $W_y\bsl\til{W}/W_x$ (which are all but finitely many), does not contribute to the formula for the orbital integral.\pp

\begin{proof}Recall $d_x-d_y>0$. Let $\Phi=\Phi(\G,\mathbf{S})$ denote the roots in the root system. For an unbounded polyhedron described above, let $\Phi^{\flat}:=\{\alpha\in\Phi\,|\,\alpha(v)\le d_x-d_y\text{ on that polyhedron}\}$. The positive span of $\alpha\in\Phi^{\flat}$ should not be the whole space, for the polyhedron has to be bounded if so. Consequently, there is a hyperplane on $X^*(\mathbf{S})\otimes\mathbb{R}$ such that all roots lies on either the boundary or one side of the hyperplane.\p

We may take such a hyperplane to have integral coefficient. In other words, we may find a cocharacter $\lambda:\mathbb{G}_m\ra\mathbf{S}$ such that $\langle\lambda,\alpha\rangle\ge 0,\;\forall\alpha\in\Phi^{\flat}$. In other words, $\alpha(x-w^{-1}y)>d_x-d_y$ for all $\alpha$ with $\langle\lambda,\alpha\rangle<0$. We may identify the cocharacter $\lambda$ as $\bar{\lambda}:\mathbb{G}_m\ra\bar{\mathbf{S}}$ just like in Lemma \ref{tech}. By (\ref{Hx}), an element in $((\LG_{w^{-1}y,d_y}+\LG_{x,d_x+})\cap\LG_{x,d_x})/\LG_{x,d_x+}$ lies on the non-negative root space of $\bar{\lambda}$. However, the proof of Theorem \ref{ellip} exactly shows that this cannot happen for elements with closed $\mathbf{L}_x$-orbit and anisotropic stabilizer.\end{proof}\pp\pp

Consider $\til{W}=\til{Z}\rtimes W$ where $\til{Z}$ is the quotient of $C_G(S)$ by the stabilizer of any point on $\mathcal{A}=\mathcal{A}(\mathbf{S},F)$ in $C_G(S)$ (this $\til{Z}$ is the group denoted by $\Lambda$ in \cite[1.2]{Tits}). $W=N_G(S)/C_G(S)$ is the finite Weyl group of $\G$ over $F$. We have $W_x$ and $W_y$ maps injectively into $W$. There is a ``derivative'' map $D:W_y\bsl\til{W}\ra W_y\bsl W$. We'll need another definition:

\begin{definition}\label{linked} Elements $w, w'\in W_y\bsl\til{W}$ with the same image $D(w)=D(w')$ in $W_y\bsl W$ is called {\bf linked} if $x-w^{-1}y$ and $x-(w')^{-1}y$ lies in the same polyhedron in $\V$. They are called {\bf strongly linked} \footnote{the reader might want to skip this notion on reading; see Remark \ref{lol}} if they are linked, and there exists a lift $z\in C_G(S)$ of the element $\bar{z}\in\til{Z}$ with $w\bar{z}=w'$ such that the action of $\Ad(z)$ on
$$\sum_{\alpha\in\Phi(\G,\mathbf{S}),\,\alpha(x-w^{-1}y)=d_x-d_y}\LG_{x,d_x:d_x+,\alpha}$$

is the identity.\p

We say elements $w,w'\in W_y\bsl\til{W}/W_x$ are linked (resp. strongly linked) if some lifts of them in $W_y\bsl\til{W}$ are linked (resp. strongly linked). We write $w\sim w'$ if they are linked and $w\approx w'$ if they are strongly linked. A {\bf linked class} (resp. {\bf strongly linked class}) is an equivalence class for $\sim$ (resp. $\approx$) in $W_y\bsl\til{W}$ or $W_y\bsl\til{W}/W_x$ whose image under $ev$ lies in a bounded polyhedron.
\end{definition}\pp\pp

Note that in the first part of the definition, the action of $\Ad(z)$ on the space is always well-defined because $w\sim w'$ implies $\alpha(x-w^{-1}y)=d_x-d_y\Leftrightarrow \alpha(x-w'^{-1}y)=d_x-d_y$ and that $\Ad(z)$ preserves $\sum_{\text{such }\alpha}\LG_{x,d_x,\alpha}$. The following lemma shows that there are at most a bounded number of strongly linked classes among a linked class.\p

\begin{lemma}\label{easy} Let $Z=\mathbf{S}(F)/\mathbf{S}(\CO_F)$, identified as a finite index subgroup of $\til{Z}$. If for $w,w'\in W_y\bsl\til{W}$ we have $w\sim w'$ and $\bar{z}\in Z$ with $w\bar{z}=w'$, then $w\approx w'$.
\end{lemma}

\begin{proof} Choose a uniformizer $\pi$ for our local field $F$. One has a natural identification $Z\cong X_*(\mathbf{S})$ by identifying the image of the discrete valuation of $F$ to $\mathbb{Z}$. This gives a lift $Z\cong X_*(\mathbf{S})\ra S=\mathbf{S}(F)$ by $\chi\mapsto\chi(\pi)$. That $w\sim w'$ implies $(\bar{z},\alpha)=0$ for all $\alpha\in\Phi(\G,\mathbf{S})$ with $\alpha(x-w^{-1}y)=d_x-d_y$. But then for all such $\alpha$, $z=\bar{z}(\pi)$ acts on $\mathfrak{u}_{\alpha}$ by $\pi^{(\bar{z},\alpha)}=1$.
\end{proof}\p

\begin{remark}\label{lol} The lemma implies that there can only be a bounded amount of strongly linked class inside a linked class. In fact, the author does not know a single example of a linked class in which there are more than one strongly linked class, and it can be proved that for a quasi-split group (i.e. when $C_G(S)$ is a torus) a linked class is always a strongly linked class. We don't know if this is true in general or not.
\end{remark}\p

Now we can finally rewrite the last equation of the last subsection. Denote by $\mathcal{C}$ the set of strongly linked classes in $W_y\bsl\til{W}/W_x$.\p

\begin{theorem} We have
\begin{equation}\label{main}
\mu_{\til{\phi}_x}(1_{\phi_y+\LG_{y,d_y+}})=\sum_{\Pi\in\mathcal{C}}\left(\sum_{w\in\Pi}q^{m_w-m_w'}\right)\cdot\left(\sum_{\dot{\gamma}\in H_{x,\Pi}}j_{\Pi,\phi_y}(\Ad(\dot{\gamma})\phi_x)\right).
\end{equation}\pp\pp

The notations are as follows: firstly $\mathbf{F}_{x,\Pi}:=\mathbf{P}_{x,w}\bsl\mathbf{L}_x$ for any $w\in\Pi$ by Lemma \ref{tech}, and $F_{x,\Pi}=\mathbf{F}_{x,\Pi}(k)$. And

$$H_{x,\Pi}=\{\dot{\gamma}\in F_{x,\Pi}\,|\,\Ad(\dot{\gamma})\phi_x\in\sum_{\alpha\in\Phi(\G,\mathbf{S}),\alpha(x-w^{-1}y)\le d_x-d_y}\LG_{x,d_x:d_x+,\alpha}\}.$$\pp

Lastly, $j_{\Pi,\phi_y}:\LG_{x,d_x:d_x+}\ra\Z_{\ge 0}$ is given by

\begin{equation}\label{aftmain}\dsp j_{\Pi,\phi_y}(a)=\#\text{Stab}_{L_{w^{-1}y}(\phi_{w^{-1}y})}\cdot\#\{\phi'\in\Ad(L_{w^{-1}y})(\phi_{w^{-1}y})\;|\;\phi'\in(\LG_{w^{-1}y,d_y+}+\LG_{x,d_x})/\LG_{w^{-1}y,d_y+}
\end{equation}
$$\text{ and its reduction in }\sum_{\alpha\in\Phi(\G,\mathbf{S}),\alpha(x-w^{-1}y)=d_x-d_y}\LG_{x,d_x:d_x+,\alpha}\text{ agrees with }a\}.$$\pp\pp

\end{theorem}

Recall the notations $\phi_{w^{-1}y}=\Ad(w^{-1})\phi_y$ and $q^{m_w-m_w'}=[\LG_{x,0+}:\LG_{x,0+}\cap\LG_{w^{-1}y}]\cdot[\LG_{x,d_x+}:\LG_{x,d_x+}\cap\LG_{w^{-1}y,d_y+}]^{-1}$. Note that for $w\in\Pi$ the definitions of $j_{\Pi,\phi_y}$ are independent of the choice of $w$, since $wz=w'\Rightarrow\phi_{w^{-1}y}=\Ad(z)\phi_{w'^{-1}y}$; this is where we use the notion of strongly linked.\p

The set $H_{x,\Pi}$ is actually the $k$-points of a subvariety $\mathbf{H}_{x,\Pi}$ of the flag variety $\mathbf{F}_{x,\Pi}$ which are called {\it Hessenberg varieties} in \cite{GKM06}. In fact, when $y$ is a hyperspecial vertex, $d_y=0$ and if we replace $1_{\phi_y+\LG_{y,d_y+}}$ in (\ref{main}) by $1_{\LG_{y,d_y}}$, then we are computing orbital integral on the indicator function of (the Lie algebra of) a hyperspecial parahoric, and what we are doing should be thought as a $p$-adic integral analogue of the work of \cite{GKM06}. In general, we should be counting points on quasi-finite covers of Hessenberg varieties. The author studies some examples in \cite[Sec. 3]{Ts15c}.\p

At this stage, the combinatorics of the set $\mathcal{C}$ of strongly linked classes can be very complicated, since the set of linked classes, or equivalently, the set of polyhedrons cut out on $\V$ is already very complicated. In the next section we'll see a situation where we simplify this combinatorics, and use the result to compute Shalika germs.\p

\section{Shalika germs}\label{Shalikagerms}

\subsection{Homogeneity of nilpotent orbits}

Let $x$, $d_x$, $\phi_x$ and $\til{\phi}_x$ any lift of $\phi_x$ be as in the previous section. Fix in this section $N$ a nilpotent element in $\LG$. The ($G$-conjugacy) orbit of $N$ will be denoted ${}^GN$. We shall assume in this section $\text{char}(F)$ is large enough so that there are only finitely many such orbits, see Appendix. We shall compute the value of the Shalika germ $\Gamma_{N}(\til{\phi}_x)$ in terms of counting points on subvarieties of flag varieties and those combinatorics same as in (\ref{main}).\p

\begin{hypothesis}\label{nil}There exists a point $y$ on the building, a real number $d_y^*$ and a coset $\phi_y^*\in\LG_{y,d_y^*}/\LG_{y,d_y^*+}$ such that $N\in\phi_y^*+\LG_{y,d_y^*+}$, and for any nilpotent orbit other than ${}^GN$ that meets $\phi_y^*+\LG_{y,d_y^*+}$, the dimension of that orbit is strictly larger than the dimension of ${}^GN$. 
\end{hypothesis}\pp\pp

By work of DeBacker \cite[Corollary 5.2.4]{De02b} the hypothesis always holds when $\text{char}(k)$ is large enough; In his construction, all other orbits that meets the coset contains ${}^GN$ in its ($p$-adic) closure.\p

We now suppose the hypothesis holds for $N$, and we fix a choice of $y$, $d_y^*$ and $\phi_y^*$ that satisfy the condition in the hypothesis. Write $\mu_{\CO}$ for the nilpotent orbital integral on a nilpotent orbit $\CO$ and $\mu_N:=\mu_{({}^GN)}$. As a consequence, $\mu_N(1_{\phi_y^*+\LG_{y,d_y^*+}})>0$, and $\mu_{\CO}(1_{\phi_y^*+\LG_{y,d_y^*+}})=0$ for a different orbit $\CO$ of the same or less dimension as that of ${}^GN$.\p

For any nilpotent element $N'$, by \cite[Proposition 1]{Po}, there is always a $1$-dimensional torus in $\G$ that acts on the line generated by $N'$ but non-trivially. Such torus must split, and thus we have a torus $j:\mathbb{G}_m\hra\G$ which acts on $N'$ by $\Ad(j(t))N'=t^{\ell_{N'}}$ for some $\ell_{N'}\in\Z-\{0\}$. When $F$ is a $p$-adic field one can integrate an $\mathfrak{sl}_2$-triple and take $\ell_{N'}=2$. We take $\ell$ to be the least common multiple for a choice of such $\ell_{N'}$ for each nilpotent orbit, so that for all nilpotent element $N'$, $t^{\ell}N'\in{}^GN'$ for any $t\in F^*$.\p

Fix also a choice of the uniformizer $\pi\in F$. We in particular have $\pi^{\ell}N'\in{}^GN'$. The assumption that $\text{char}(F)$ is very good for $\G$ guarantees that the centralizer of $N'$ in $\LG$ is the Lie algebra of the centralizer of $N'$ in $G$ \cite[Theorem 1.2]{BMRT} (this implies \cite[Proposition 6.7]{Bo} that the nilpotent orbit is smooth with the expected tangent space).\p

Recall by \cite[Proposition 4.1]{AR00} we have a $G$-invariant non-degenerate symmetric bilinear form $\kappa(\cdot,\cdot)$ on $\LG$. One can then consider the symplectic form on ${}^GN'$ (as an $F$-analytic manifold over $F$), for which this form at the point $N'\in{}^GN'$ is given by $\langle X,Y\rangle_{N'}=\kappa(\Ad(N')X,Y)$ for $X,Y\in\LG/\LG_{N'}$, where $\LG_{N'}$ is the Lie algebra of the centralizer of $N'$.\p

The top exterior power of this symplectic form give rises to an $G$-invariant measure on ${}^GN'$. Obviously $\langle\cdot,\cdot\rangle_{t^\ell N'}=t^\ell\langle\cdot,\cdot\rangle_{N'}$ for $t\in F^*$. This implies that the measure on the nilpotent orbit satisfies the property that for any subset $S\subset\LG$, the measure of ${}^GN'\cap\pi^\ell S$ is $q^{-\frac{\dim{}^GN'}{2}\cdot\ell}$ times the measure of ${}^GN'\cap S$.\p

For any function $f\in C_c^{\infty}(\LG)$, $t\in F^*$, denote by $f_t$ the function $f_t(X)=f(t^{-1}X)$. Then for any nilpotent orbit $\CO$
$$\mu_{\CO}(f_{\pi^{\ell n}})=q^{-\frac{\dim\CO}{2}\cdot\ell n}\mu_{\CO}(f),\;\forall n\in\Z.$$\pp

The theorem of Shalika germs \cite{Sh72} states that under certain assumptions on $\text{char}(F)$ (see Appendix), there exists smooth functions $\Gamma_{\CO}$ on $\LG^{rs}$ the dense open subset of regular semisimple elements of $\LG$, such that for any $X\in\LG^{rs}$, there exists a lattice $\Lambda_X\subset\LG$ such that for any compactly supported function $f\in C_c^{\infty}(\LG/\Lambda_X)$ invariant under translation by $\Lambda_X$, we have
$$\mu_X(f)=\sum_{\CO}\Gamma_{\CO}(X)\mu_{\CO}(f).$$

Here $\CO$ runs over nilpotent orbits, and $\mu_X$ and $\mu_{\CO}$ are orbital integrals for the orbit of $X$ and $\CO$, respectively, with suitable normalization for their Haar measures. For nilpotent element $N$ we also write $\Gamma_N:=\Gamma_{({}^GN)}$.\p

Let $y$, $d_y^*$ and $\phi_y^*$ be as in Hypothesis \ref{nil}. Let $n$ be a positive integer. In the rest of this section let $d_y=d_y^*+\ell n$ and $\phi_y=\pi^{-\ell n}\phi_y^*$ depend on $n$. We consider $\mu_{\til{\phi}_x}(1_{\phi_y+\LG_{y,d_y+}})$. Since multiplying by $\pi^{-\ell n}$ preserves all nilpotent orbits, the new pair of $d_y$ and $\phi_y$ still satisfy Hypothesis \ref{nil}. Together with the theorem of Shalika germ, we have

\begin{equation}\label{Shalika}\dsp\mu_{\til{\phi_x}}(1_{\phi_y+\LG_{y,d_y+}})=\sum_{\CO}\Gamma_{\CO}(\til{\phi}_x)\mu_{\CO}(1_{\phi_y+\LG_{y,d_y+}})=\sum_{\CO}q^{\frac{\dim\CO}{2}\cdot \ell n}\Gamma_{\CO}(\til{\phi}_x)\mu_{\CO}(1_{\phi_y+\LG_{y,d_y+}}).
\end{equation}\pp\pp

When $\CO={}^GN$,  $\frac{\dim\CO}{2}=\frac{\dim({}^GN)}{2}$. Hypothesis \ref{nil} says that for all other orbits $\CO'$ that appears in this sum, $\frac{\dim\CO'}{2}>\frac{\dim({}^GN)}{2}$. We shall compare (\ref{Shalika}) with (\ref{main}) in the next section. 

\subsection{Relatedness and strong-relatedness}

Let $n'$ be a positive integer and $d_y'=d_y^*+\ell n'$, $\phi_y'=\pi^{-\ell n'}\phi_y^*$.\p

\begin{definition}\label{relate}
Let $\W_{d_y}^{\sim}$ denote the set of linked classes in $W_y\bsl\til{W}$ (while the definition depends on $d_y$ and $n$) and $\W_{d_y'}^{\sim}$ instead be the set of linked classes defined with $d_y$ replaced by $d_y'$. We say $\Pi\in\W_{d_y}^{\sim}$ and $\Pi'\in\W_{d_y'}^{\sim}$ are {\bf related} if any (hence all) $w\in\Pi$, $w'\in\Pi'$ we have $D(w)=D(w')$, and for any $\alpha\in\Phi(G,S)$ one has $\alpha(x-w^{-1}y)$ and $\alpha(x-w'^{-1}y)$ have the same sign and $\alpha(x-w^{-1}y)-(d_x-d_y)$ and $\alpha(x-w'^{-1}y)-(d_x-d_y')$ have the same sign.\p

Similarly, let $\W_{d_y}^{\cong}$ and $\W_{d_y'}^{\cong}$ denote the set of strongly linked classes in $W_y\bsl\til{W}$, defined with $d_y$ and $d_y'$, respectively. We say $\Pi\in\W_{d_y}^{\cong}$ and $\Pi'\in\W_{d_y'}^{\cong}$ are {\bf strongly related} if they belong to related linked classes, and there exists $\bar{z}\in\til{Z}$ and its lift $z\in C_G(S)$ such that for some $w\in\Pi$, $w'\in\Pi'$ one has $w\bar{z}=w'$ and that $\Ad(z)$ acts on
$$\bigoplus_{\alpha\in\Phi(G,S),\,\alpha(x-w^{-1}y)=d_x-d_y}\LG_{x,d_x:d_x+,\alpha}=\bigoplus_{\alpha\in\Phi(G,S),\,\alpha(x-w'^{-1}y)=d_x-d_y'}\LG_{x,d_x:d_x+,\alpha}$$

by $\pi^{kn}$. Lastly, we write $\mathcal{C}_{d_y}^{\sim}:=\W_{d_y}^{\sim}/W_x$ and $\mathcal{C}_{d_y}^{\cong}:=\W_{d_y}^{\cong}/W_x$, the same for $d_y'$ and induce the notion of relatedness and strong relatedness on them. That is, elements in $\mathcal{C}_{d_y}^{\sim}$ and $\mathcal{C}_{d_y'}^{\sim}$ (resp. $\mathcal{C}_{d_y}^{\cong}$ and $\mathcal{C}_{d_y'}^{\cong}$) are related (resp. strongly related) if some of their preimages are.\p
\end{definition}\pp

Exactly the same argument as in Lemma \ref{easy} shows that if $w\bar{z}=w'$ for related $w\in\Pi$, $w'\in\Pi'$ and $\bar{z}\in Z$, then $\Pi$ and $\Pi'$ are strongly related. This notion of strong relatedness is for the reason that, just like the notion for strongly linkedness, the geometrical sum (i.e. the sum in the last round bracket) in (\ref{main}) depends only up to strong relatedness.\p

Recall in the paragraph before Lemma \ref{finite} we took $\V=X_*(\mathbf{S})\otimes\mathbb{R}$. Consider this time affine hyperplanes $\alpha(v)=0$ and $\alpha(v)=2$ for $\alpha\in\Phi(\G,\mathbf{S})$. The affine hyperplanes cut $\V$ into numerous polyhedrons (of different dimensions). Call $\CP$ the set of these polyhedrons. The group $W_x$ acts on the space $\V$ by fixing the origin. Denote by $\bar{\CP}$ the set of $W_x$-orbits in $\CP$. By the very definition of relatedness we have injective maps $ev:\W_{d_y}^{\sim}\ra\CP$ and $ev:\CC_{d_y}^{\sim}\ra\bar{P}$ by $w\mapsto\left[\frac{2}{d_x-d_y'}(x-w^{-1}y)\right]$.\p

By Lemma \ref{finite}, only those (strongly) linked class that maps to bounded polyhedrons matter. Denote by $\W_{d_y}^{\sim,\text{bd}}\subset\W_{d_y}^{\sim}$ the subset of those which map to bounded polyhedrons and $\CC_{d_y}^{\sim,\text{bd}}\subset\CC_{d_y}^{\sim}$ those which map to orbits of bounded polyhedrons. As we have natural map $\W_{d_y}^{\cong}\ra\W_{d_y}^{\sim}$ and $\CC_{d_y}^{\cong}\ra\CC_{d_y}^{\sim}$, we can define $\W_{d_y}^{\cong,\text{bd}}$ and $\CC_{d_y}^{\cong,\text{bd}}$, and similarly for $d_y'$.\p

Recall that in (\ref{aftmain}), we have
\begin{equation}\label{main2}\mu_{\til{\phi_x}}(1_{\phi_y+\LG_{y,d_y+}})=\sum_{\Pi\in\CC_{d_y}^{\cong,\text{bd}}}\left(\sum_{w\in\Pi}q^{m_w-m_w'}\right)J(\Pi).
\end{equation}

where $J(\Pi)$, being the second inner sum in the RHS of (\ref{main}), is some ``geometric number'' that only depends on $\Pi$ up to strongly relatedness. We now have\p

\begin{theorem}\label{combina} After possibly replacing $\ell$ by a multiple of it, we have\pp

(1) There is an integer $N$ such that for all $n,n'\ge N$, we have bijections $\mathcal{C}_{d_y}^{\sim,\text{bd}}\xra{\sim}\mathcal{C}_{d_y'}^{\sim,\text{bd}}$ and $\mathcal{C}_{d_y}^{\cong,\text{bd}}\xra{\sim}\mathcal{C}_{d_y'}^{\cong,\text{bd}}$ given by relatedness and strong relatedness.\pp
(2) Let $\mathcal{C}^{\sim,\text{bd}}$ be the set identified to all $\mathcal{C}_{d_y}^{\sim\text{bd}}$ and similarly by $\mathcal{C}^{\cong\text{bd}}$ the set identified to all $\mathcal{C}_{d_y}^{\cong}$ as in (1). Then we have for any strongly related class $\Pi\in\mathcal{C}^{\cong,\text{bd}}$, the sum $\sum_{w\in\Pi}q^{m_w-m_w'}$ is a polynomial $p_{\Pi}(q^n,n)$ in $q^n$ and $n$.
\end{theorem}\pp

\begin{proof} For (1), since for any $\sigma\in W_x$, $w\in\W_{d_y}^{\sim,\text{bd}}$ and $w'\in\W_{d_y'}^{\sim,\text{bd}}$ are related if and only if $w\sigma$ and $w'\sigma$ are related (and the same for strongly related), it suffices to find a multiple of $\ell$ so that $\W_{d_y}^{\sim,\text{bd}}\xra{\sim}\W_{d_y'}^{\sim,\text{bd}}$ and $\W_{d_y}^{\cong,\text{bd}}\xra{\sim}\W_{d_y'}^{\cong,\text{bd}}$ under relatedness and strong relatedness. We need an extension of Lemma \ref{easy}, whose proof is the same:\pp

\begin{lemma}\label{easy2} Let $Z\subset\til{Z}$ be as in Lemma \ref{easy}. If for $\Pi\in\W_{d_y}^{\cong}$, $\Pi'\in\W_{d_y'}^{\cong}$, $\Pi$ and $\Pi'$ are related and for some (hence all) $w\in\Pi$, $w'\in\Pi'$ there exists $\bar{z}\in Z$ with $w\bar{z}=w'$, then $\Pi$ and $\Pi'$ are strongly related.
\end{lemma}\p

Let $\W_{d_y}^{\equiv}$ be a further partition of $\W_{d_y}^{\sim}$ by the property in Lemma \ref{easy}. The above lemma stated that $\W_{d_y}^{\equiv}$ is also a partition of $\W_{d_y}^{\cong}$. Let's invent a new term that $\Pi\in\W_{d_y}^{\equiv}$ and $\Pi'\in\W_{d_y'}^{\equiv}$ are called equivalently related if for some $w\in\Pi$, $w'\in\Pi'$, the property in the lemma holds. Then it suffices to prove the bijection $\W_{d_y}^{\equiv,\text{bd}}\xra{\sim}\W_{d_y'}^{\equiv,\text{bd}}$ under equivalent relatedness.\p

Note that $w\mapsto x-w^{-1}y$ gives an injection from $W_y\bsl\til{W}$ to $\mathcal{A}$. For elements $w,w'\in\Pi$ for some $\Pi\in\W_{d_y}^{\equiv,\text{bd}}$, we have $x-w^{-1}y$ and $x-w'^{-1}y$ satisfy that $\alpha(x-w^{-1}y)$ and $\alpha(x-w'^{-1}y)$ shares the same sign, that $\alpha(x-w^{-1}y)-(d_x-d_y)$ and $\alpha(x-w'^{-1}y)-(d_x-d_y')$ shares the same sign, and that $x-w^{-1}y$ and $x-w'^{-1}y$ differs by a translation from $Z$, where $Z$ acts on $\mathcal{A}$ by translation through a lattice. The assertion then follows from the following combinatorial lemma:

\begin{lemma}\label{comb1}We use notations independent from the rest of the article. Consider lattice $\Z^m\subset\mathbb{R}^m$ and a finite set of linear equations and inequalities of the form $a_1x_1+...+a_nx_m=c(n)$ and $a_1x_1+...+a_mx_m<c(n)$, where $a_i$ are rational numbers and $c(n)$ is a linear function in integer variable $n$ with rational coefficients. Suppose the equations and inequalities satisfy that the region of solutions in $\mathbb{R}^m$ is bounded. Let $q$ be a formal variable and $\lambda(\mathbf{x})$ be a linear function on $\Z^m$ with integral coefficients. Let $S_{\lambda}(n)\in\Z[q,q^{-1}]$ be the sum of $q^{\lambda(\mathbf{x})}$ over all $\mathbf{x}\in\Z^m$ in the region, i.e. that satisfy the equations and inequalities. Then there exists positive integers $c$ and $N$ such that for all $n>N$, $S_{\lambda}(cn)$ is a Laurent polynomial in $q^n$ and $n$.
\end{lemma}\pp

One takes the conditions regarding the signs of $\alpha(x-w^{-1}y)$ and $\alpha(x-w^{-1}y)-(d_x-d_y)$ to be the equations and the inequalities. The lemma implies that the number of lattice points inside such a set of equations and inequalities is either always non-zero or always zero (one sees this by subtituting $q=1$). Hence the result for (1).\p

For (2) of the theorem, the lemma almost provide the result except that we are summing over $w$ in subsets of $W_y\bsl\til{W}/W_x$ instead of $W_y\bsl\til{W}$ (as in the setting of Definition \ref{relate}, Lemma \ref{easy2} and the above lemma. This difference can be resolved by adding extra linear equations to ``mark'' the fixed points of various subgroups of $W_x$, so that in each subset $\Pi$ of $W_y\bsl\til{W}$ considered, $W_x$ acts on the right with the order of the stabilizer being constant, and the sum over $\Pi\subset W_y\bsl\til{W}/W_x$ becomes just the sum over its preimage in $W_y\bsl\til{W}$ divided by this order. This finishes the proof of Theorem \ref{combina}.
\end{proof}\p

Let's now fix $\Pi\in\mathcal{C}^{\cong,\text{bd}}$. We analyze the polynomial $p_{\Pi}(q^n,n)$ in terms of $ev(\Pi)\in\bar{\CP}$. Recall that in (\ref{main}) $q^{m_w-m_w'}$ was given by 

$$q^{m_w-m_w'}=[\LG_{x,0+}:\LG_{x,0+}\cap\LG_{w^{-1}y}]\cdot[\LG_{x,d_x+}:\LG_{x,d_x+}\cap\LG_{w^{-1}y,d_y+}]^{-1}$$
$$=\prod_{\alpha\in\Phi(\G,\mathbf{S})}\frac{|\LG_{x,0+,\alpha}|\cdot|\LG_{x,d_x+,\alpha}\cap\LG_{w^{-1}y,d_y+,\alpha}|}{|\LG_{x,0+,\alpha}\cap\LG_{w^{-1}y,0,\alpha}|\cdot|\LG_{x,d_x+,\alpha}|},\;\;\;\;\;$$\pp

where $|\cdot|$ was used to denote any Haar measure on the $F$-vector space $\LG_{\alpha}$. We have following easy estimate\pp\pp

\begin{lemma}\label{esti}For any point $z$ in the Bruhat-Tits building, any real numbers $d_1,d_2$ and any $\alpha\in\Phi(\G,\mathbf{S})$. Let $\rho_{\alpha}=\dim_F(\LG_{\alpha})$. Suppose $[\LG_{z,d_1,\alpha}:\LG_{z,d_1,\alpha}\cap\LG_{z,d_2,\alpha}]=q^m$. Then
$$|m-\rho_{\alpha}\cdot\max(d_2-d_1,0)|<\rho_{\alpha},$$
\end{lemma}\p

Applying the estimate, we have
$$m_w-m_w'=\sum_{\alpha\in\Phi(\G,\mathbf{S})}\rho_{\alpha}\cdot\left(\max\{\alpha(x-w^{-1}y),0\}-\max\{\alpha(x-w^{-1}y)-(d_x-d_y),0\}\right)+\epsilon,$$\pp\pp

with the error term $|\epsilon|<\sum_{\alpha}\rho_{\alpha}=\dim\G-\dim C_{\G}(\mathbf{S})$. Recall that $d_y=d_y^*-\ell n$ depends on $n$. Suppose $n$ is large and $w\in W_y\bsl\til{W}$ is such that $\frac{2}{d_x-d_y'}(x-w^{-1}y)$ is very close to a vertex $v$, i.e. $0$-dimensional polyhedron in $\CP$, then one has $m_w-m_w'\sim\frac{\ell n}{2}\cdot\tau(v)$, where $\tau(v):=\sum_{\alpha\in\Phi(\G,\mathbf{S})}\rho_{\alpha}\tau_{\alpha}(v)$, with
$$\tau_{\alpha}(u)=\left\{\begin{matrix}
2&\text{if }\alpha(v)\ge 2.\\
\alpha(v)&\;\;\;\;\;\,\text{ if }0\le\alpha(v)\le 2.\\
0&\text{otherwise. }
\end{matrix}\right.$$

More precisely, we have

\begin{theorem}\label{degree}For each monomial in the polynomial $p_{\Pi}(q^n,n)$, the degree in $q^n$ is equal to $\frac{\ell\cdot\tau(v)}{2}$, for some vertex $v\in\bar{\CP}$ contained in the closure of $ev(\Pi)$ (see the definition of $ev(\Pi)$ before Theorem \ref{combina}.)
\end{theorem}\pp


\begin{proof}To prove the theorem, one simply use the following extension of the combinatorial lemma \ref{comb1}:\pp

\begin{lemma}\label{comb2}We use notations in Lemma \ref{comb1} and independent from all other parts of the article. Suppose $d$ is the degree in $q^n$ of a term in $S_{\lambda}(cn)$ in Lemma \ref{comb1}. Then possibly after replacing $c$ by a multiple of it, there exist vectors $\mathbf{v},\mathbf{w}\in\mathbb{R}^m$ and $\epsilon\in\mathbb{R}$, such that $d\cdot cn=\lambda(\mathbf{v}+cn\mathbf{w})-\epsilon$, and satisfy:\pp
(1) By changing some of the inequalities into equations (by turning the $<$ sign into $=$ sign) and throwing away other inequalities, one has that for every $n$, $\lambda(\mathbf{v}+n\mathbf{w})$ is the unique solution in $\mathbb{R}^m$ to the resulting set of equations.\pp
(2) There exists an $N'$ such that for all $n>N'$ that are divided by $c$, $\lambda(\mathbf{v}+n\mathbf{w})$ is contained in the closure of the region cut out by the equations and inequalities.
\end{lemma}

The theorem is then proved by using Lemma \ref{esti} and substituting $\lambda(x-w^{-1}y)=m_w-m_w'$ as when we used Lemma \ref{comb1}.\end{proof}\p

It's now time to compare (\ref{Shalika}) with (\ref{main2}). They are

$$\mu_{\til{\phi_x}}(1_{\phi_y+\LG_{y,d_y+}})=\sum_{\CO}q^{\frac{\dim\CO}{2}\cdot \ell n}\Gamma_{\CO}(\til{\phi}_x)\mu_{\CO}(1_{\phi_y+\LG_{y,d_y+}}).$$
$$\mu_{\til{\phi_x}}(1_{\phi_y+\LG_{y,d_y+}})=\sum_{\Pi\in\CC_{d_y}^{\cong,\text{bd}}}p_{\Pi}(q^n,n)J(\Pi).$$\p

Recall our goal is to compute $\Gamma_N(\til{\phi}_x)$, and Hypothesis \ref{nil} says $\mu_{\CO}(1_{\phi_y+\LG_{y,d_y+}})\not=0$ only when $\CO={}^GN$ or $\dim\CO>\dim{}^GN$. Comparing two equations gives

\begin{theorem}\label{mainthm}Let $\phi_x\in\LG_{x,d_x:d_x+}$ be with closed orbit and anisotropic stabilizer under $\mathbf{L}_x$-action, $\til{\phi}_x\in\LG_{x,d_x}$ any lift of $\phi_x$, and $y$, $d_y$, $\phi_y$ and $N$ as in Hypothesis \ref{nil}. Then

$$\Gamma_N(\til{\phi}_x)=(\mu_N(1_{\phi_y+\LG_{y,d_y+}}))^{-1}\cdot\sum_{\Pi\in\CC^{\cong,\text{bd}}}c_{\frac{\ell\cdot\dim\CO}{2}}(p_{\Pi})J(\Pi),$$\pp

where $c_{\frac{\ell\cdot\dim\CO}{2}}(p_{\Pi})$ denote the coefficient of the $q^{\frac{\ell\cdot\dim\CO}{2}n}$-term in $p_{\Pi}(q^n,n)$, and $J(\Pi)$ is the second inner sum in the RHS of (\ref{main}).
\end{theorem}

Logically Theorem \ref{mainthm} is independent from Theorem \ref{degree}. However, Theorem \ref{degree} gives a necessary condition for $c_{\frac{\ell\cdot\dim\CO}{2}}(p_{\Pi})$ to be non-zero; it can be non-zero only when a vertex $v$ on the boundary of $ev(\Pi)$ has $\tau(v)=\dim\CO$. In the next subsection we give another strong result along this line.\p

\subsection{A vanishing result}

This subsection is devoted to the following result.

\begin{theorem}\label{vanish}Let $\Pi\in\CC^{\cong,\text{bd}}$. Suppose there is a vertex $v$ on the boundary of $ev(\Pi)$ such that $\tau(v)<\dim\CO$. Then $J(\Pi)=0$.
\end{theorem}

For our convenience, we copy from (\ref{main}), (\ref{aftmain}) that

$$J(\Pi)=\sum_{\dot{\gamma}\in F_{x,\Pi}}j_{\Pi,\phi_y}(\Ad(\dot{\gamma})\phi_x)$$\pp

where
$$j_{\Pi,\phi_y}(a)=\#\text{Stab}_{L_{w^{-1}y}(\phi_{w^{-1}y})}\cdot\text{id}\left(a\in\sum_{\alpha\in\Phi(\G,\mathbf{S}),\alpha(x-w^{-1}y)\le d_x-d_y}\LG_{x,d_x:d_x+,\alpha}\right)\cdot$$
$$\#\{\phi'\in\Ad(L_{w^{-1}y})(\phi_{w^{-1}y})\;|\;\phi'\in(\LG_{w^{-1}y,d_y+}+\LG_{x,d_x})/\LG_{w^{-1}y,d_y+}\text{ and its reduction in }$$
$$\s\s\s\s\s\sum_{\alpha\in\Phi(\G,\mathbf{S}),\alpha(x-w^{-1}y)=d_x-d_y}\LG_{x,d_x:d_x+,\alpha}\text{ agrees with }a\}.$$\p

In particular, we see that in order to have $J(\Pi)\not=0$, we need to have

$$\phi'\in((\LG_{w^{-1}y,d_y+}+\LG_{x,d_x})\cap\LG_{w^{-1}y,d_y})/\LG_{w^{-1}y,d_y+}=\sum_{\alpha\in\Phi(\G,\mathbf{S}),\alpha(x-w^{-1}y)\ge d_x-d_y}\LG_{w^{-1}y,d_y:d_y+,\alpha}.$$

where $\phi'$ is some element in $\Ad(L_{w^{-1}y})(\phi_{w^{-1}y})$. We can then take a lift 
$N\in\LG_{w^{-1}y,d_y}$ of $\phi'$ such that

\begin{equation}\label{grade}
N\in\sum_{\alpha\in\Phi(\G,\mathbf{S}),\alpha(x-w^{-1}y)\ge d_x-d_y}\LG_{\alpha}.
\end{equation}\p

Since $d_x-d_y>0$, $N$ is nilpotent. $\phi'\in\Ad(L_{w^{-1}y})(\phi_{w^{-1}y})$ implies that for some $g\in G_y$, $\pi^{-\ell n}\Ad(g)\Ad(w)N\in\phi_y^*+\LG_{y,d_y^*+}$. By Hypothesis \ref{nil}, the dimension of the orbit of $\pi^{-\ell n}\Ad(g)\Ad(w)N$ is larger than or equal to $\dim\CO$. By definition of $\ell$, $\pi^{-\ell n}\Ad(g)\Ad(w)N$ and $N$ have the same orbit. To prove the theorem, it suffices to show that for every vertex $v$ on the boundary of $ev(\Pi)$ we have $\tau(v)\ge\dim{}^GN$.\p

Since the polyhedron $ev(\Pi)$ (in fact a $W_x$-orbit of polyhedrons, while the $W_x$-action does not matter here) is cut out by several affine walls of the form $\alpha(v)=0$ or $\alpha(v)=2$, and $w$ and $\Pi$ is related by $ev(w)=\frac{2}{d_x-d_y}(x-w^{-1}y)\in ev(\Pi)$, we see that for every point in the closure of the polyhedron, in particular for every vertex $v$ on the boundary of $ev(\Pi)$.\p

While $v$ is a $W_x$-orbit of points on $\CV=X_*(\mathbf{S})\otimes\mathbb{Q}$, we can think of it as a point on $\CV$ and thus a rational cocharacter, giving a $\mathbb{Q}$-grading on $\LG$. We write $\LG=\bigoplus_{i\in\mathbb{Q}}\LG(i)$. By (\ref{grade}), $N\in\LG(\ge2)$. We have\p

\begin{lemma}Let $\mathfrak{H}=\bigoplus\mathfrak{H}(i)$ be an $\mathbb{R}$-grading of a Lie algebra $\mathfrak{H}$ (over an arbitrary field) with $X\in\mathfrak{H}(\ge2):=\sum_{i\ge2}\mathfrak{H}(i)$. Then for every $r\in\mathbb{R}$,

$$\dim\{Y\in\mathfrak{H}\,|\,[Y,X]=0\}\ge\dim\sum_{r\le i<r+2}\mathfrak{H}(i).$$
\end{lemma}\pp

\begin{proof}The LHS is greater than or equal to the dimension of the kernel of $\text{ad}(X)$ on $\mathfrak{H}(\ge r)$, and the projection of $\text{ad}(X)|_{\mathfrak{H}(\ge r)}$ to the space on the RHS is the zero map.
\end{proof}\p

On the other hand, we have by definition
$$2\tau(v)=\tau(v)+\tau(-v)=\sum_{i\in\mathbb{Q}}\min(|i|,2)\cdot\dim\LG(i)=2\dim\LG-\int_{-2}^{0}\dim\sum_{r\le i<r+2}\LG(i)dr.$$\p

With the above lemma for $\mathfrak{H}=\LG$, this says $\tau(v)\ge\dim{}^GN$ as asserted.\p

\appendix
\section{Assumptions on $p$}\label{app}

In this appendix, we summary assumptions on either $\text{char}(F)$ or $\text{char}(k)$ at various places of this article.\p

To begin with, $\text{char}(F)$ is assumed to be very good for a lot of purposes, including even the well-definedness of orbital integrals \cite[III.3.27]{SS70}. Recall we say a prime $p$ is good ($0$ is always very good) for $\G$ if $p$ does not divide the coefficient of any simple root in the highest roots, for the absolute root system of any simple factor of (the isogeny type) of $\G$. $p$ is very good for $\G$ if $p$ is good and it does not divide $n+1$ whenever $\G$ (up to isogeny) has a simple factor of type $A_{n}$.\p

For the proof of Theorem \ref{ellip} and Claim \ref{ellipsurj} we require that $\G$ is tamely ramified, that is, $\G$ splits over a tamely ramified extension. We also assume that the rational coordinate of our point $x$ on the Bruhat-Tits building has denominator coprime to $\text{char}(k)$. However, that $\LG_{x,d_x}$ has semisimple $\mathbf{L}_x$-orbit implies that $x$ must be a ($0$-dimensional) intersection of walls defined by affine roots. Both tameness condition (on $\G$ and on $x$) will thus be automatic if $\text{char}(k)$ is larger than a certain number determined by the absolute root system of $\G$.\p

In the proof of Claim \ref{ellipsurj} we need $\text{char}(k)$ to be very good to identify $\LG$ with its dual in a way compatible with the Moy-Prasad filtration. It might be that in some cases the method in \cite[Sec. 5]{GKM06} will allow us to drop this assumption.\p

In Section \ref{Shalikagerms}, if $\text{char}(F)>0$ we need to verify the assumptions (U1)$\sim$(U4) in \cite{Sh72}. (U1) and (U2) are implied by that $\text{char}(F)$ is very good \cite[Theorem 1.2]{BMRT}, \cite[Proposition 6.7]{Bo} and \cite[III.3.27]{SS70}. (U3) is proved by Ranga Rao's \cite{Ra72}, which was stated for $\text{char}(F)=0$ but in fact only uses (U1), (U2) and \cite[III.4.14]{SS70}. The last was proved under the assumption $\text{char}(F)\ge 4m+3$ where $m$ is the sum of coefficients for the highest root in the absolute root system (for each simple factor). For (U4), the number of geometric (stable) nilpotent orbits is always finite \cite{HS}. For the number of orbits to be finite, we need the first Galois cohomology of the centralizer of any nilpotent element to be finite. If $\text{char}(F)$ is larger than $\text{rk}_{\bar{F}}(\G)+1$ where $\text{rk}_{\bar{F}}(\G)$ denotes the absolute rank of $\G$, then the Levi factor of the centralizer is always tamely ramified and consequently the $H^1$ is finite.\p

Lastly, for the computation of Shalika germs of a chosen nilpotent orbit we uses Hypothesis \ref{nil}. The hypothesis is valid if DeBacker's result \cite{De02b} is available, which is good when $\text{char}(k)$ is large enough; we refer the reader to \cite[Sec. 4.2]{De02b}. For applications in not-so-large residue characteristic one can usually check Hypothesis \ref{nil} directly.\p

\let\oldthebibliography=\thebibliography
\let\endoldthebibliography=\endthebibliography
\renewenvironment{thebibliography}[1]{
  \begin{oldthebibliography}{#1}
    \setlength{\itemsep}{2pt}
  }{
    \end{oldthebibliography}
  }
\bibliographystyle{acm}
\bibliography{biblio}

\end{document}